\documentclass[12pt,twoside]{article}
\usepackage{amsmath}
\usepackage{amsfonts}
\usepackage{amsthm}
\usepackage{amssymb}
\usepackage{graphicx}
\usepackage{multicol}

\begin{document}
\title{{\normalsize{\bf The crossing numbers of knots and links via 4D topology}}}
\author{{\footnotesize Akio KAWAUCHI}\\
{\footnotesize{\it Osaka Central Advanced Mathematical Institute, 
Osaka Metropolitan University}}\\
{\footnotesize{\it Sugimoto, Sumiyoshi-ku, Osaka 558-8585, Japan}}\\
{\footnotesize{\it kawauchi@mu.ac.jp}}}
\date\, 
\maketitle
\vspace{0.25in}
\baselineskip=10pt
\newtheorem{Theorem}{Theorem}[section]
\newtheorem{Conjecture}[Theorem]{Conjecture}
\newtheorem{Lemma}[Theorem]{Lemma}
\newtheorem{Sublemma}[Theorem]{Sublemma}
\newtheorem{Proposition}[Theorem]{Proposition}
\newtheorem{Corollary}[Theorem]{Corollary}
\newtheorem{Claim}[Theorem]{Claim}
\newtheorem{Definition}[Theorem]{Definition}
\newtheorem{Example}[Theorem]{Example}

\begin{abstract} The diagrams of  links (including knots) are characterized in terms of circular rigid disk-chord diagrams of  their spun ribbon torus-links in the 4-sphere. 
As a result, the crossing numbers  of  links are equal to the chord indexes of their spun ribbon torus-links.
By using this result,  additivity on the crossing numbers of links under connected sums 
can be shown.

\phantom{x}



\end{abstract} 

\baselineskip=15pt

\bigskip

\noindent{\bf 1. Introduction}

\begin{figure}[hbtp]
\begin{center}
\includegraphics[width=14.5cm, height=8cm]{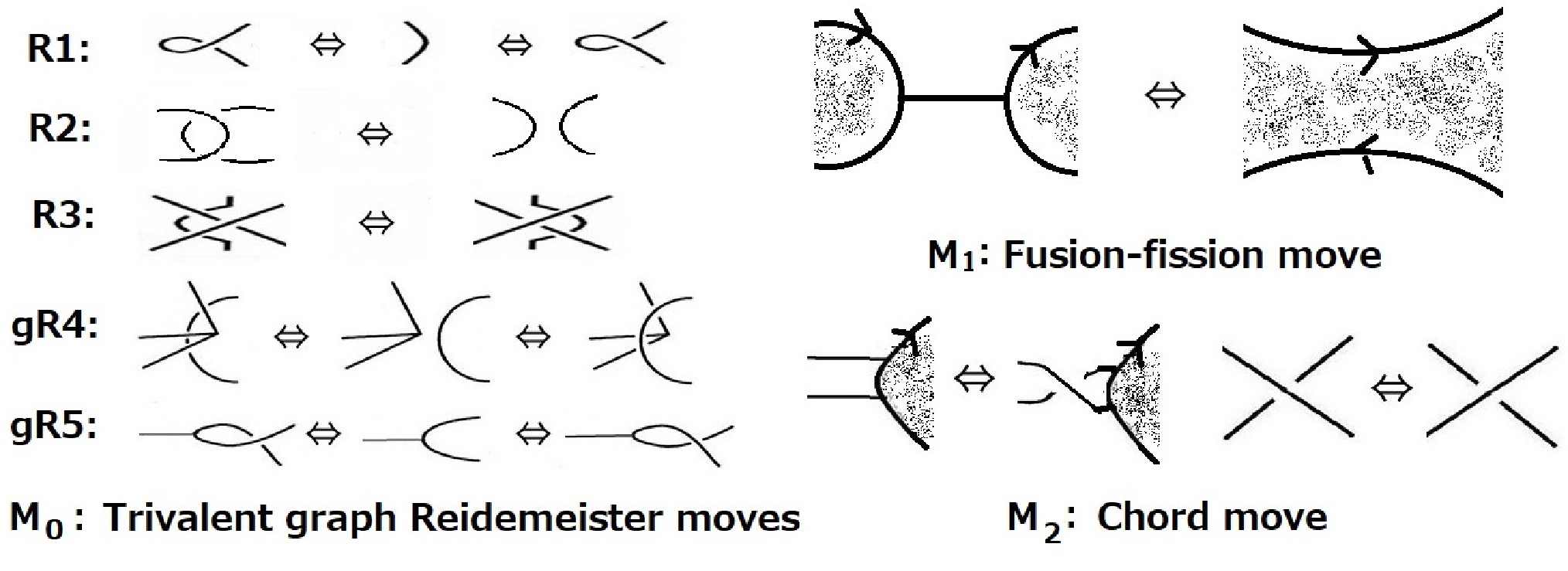}
\end{center}
\caption{The moves for equivalence of loop-chord diagrams}
\label{fig:equiv}
\end{figure}

A  {\it diagram} $D$ on the 2-sphere $S^2$ 
is an immersed oriented loop system on the 2-sphere $S^2$ whose singularities are 
only transverse double points with upper-lower relations. 
When the 2-sphere $S^2$ is fixed in the 3-sphere $S^3$, the diagram $D$  represents 
a unique link $k$ in $S^3$, obtained by pushing the interior of an upper arc of every double point of $D$ to a positive direction of $S^2$ in $S^3$. Every link in $S^3$ up to equivalences 
correspond to a diagram on $S^2$ up to Reidemiester moves, \cite{[1]}. 
A surface-link  in the 4-sphere $S^4$ is 
a (possibly disconnected) closed oriented surface $F$  smoothly embedded in  $S^4$. 
Let $r=r(F)$ be the number of components of $F$. When $r=1$,   it is called a 
surface-knot. 
The surface-link $F$ is an $S^2$-{\it link} or a {\it torus-link} if all the components of $F$ are 2-spheres or tori, respectively. 
A {\it handled-sphere system} is a pair $(O, h)$ where  $O$ is a trivial
$S^2$-link  and $h$ is a 1-handle system on $O$, smoothly embedded in $S^4$.
A {\it ribbon} surface-link $F$ is a surface-link in $S^4$ constructed from 
a handled-sphere system $(O, h)$ by surgery along $h$. 
The handled-sphere system $(O, h)$ is uniquely constructed from a 
{\it chorded-sphere system} $(O, \alpha)$ in $S^4$ where $\alpha$ is a core arc system of  the 1-handle $h$, called a {\it chord system}, \cite{[2]}.
The chorded-sphere system $(O, \alpha)$ is uniquely constructed from a 
{\it loop-chord graph} $(o, \alpha)$ located in a fixed  standard 3-sphere $S^3$ in $S^4$ 
with $o$ a trivial link, called a {\it based loop system}, 
by using the Horibe-Yanagawa's lemma, \cite{[3]}. 
A  regular diagram on  $S^2$ of a chord graph $(o, \alpha)$ in $S^3$   is a {\it loop-chord diagram}  $C(o,\alpha)$ on $S^2$, \cite{[4],[5]}. 
By convention, the based loop system $o$ on  $S^2$ 
is oriented in a clockwise orientation and situated  to bound a disjoint oriented disk system $d$ in $S^2$, called the {\it based disk system}. 
Under this convention,  the loop-chord diagram  $C(o,\alpha)$ is also considered 
as a {\it disk-chord diagram} and denoted by $C(d,\alpha)$. 
Every ribbon surface-link up to equivalences is characterized by a loop-chord diagram 
$C(o,\alpha)$ on $S^2$ up to {\it equivalences}, namely up to finite numbers of the 
following moves:\, 
Trivalent graph Reidemeister move $M_0$, Fusion-fission move $M_1$ and Chord move $M_2$, illustrated in Fig.~\ref{fig:equiv}, \cite{[4], [5],  [6], [7]}. 
Two disk-chord diagrams $C(d,\alpha)$ and $C(d',\alpha')$ are {\it homotopic} if 
a neighborhood of $d$ in $C(d,\alpha)$ is isotopic to  a neighborhood of $d'$ in 
$C(d',\alpha')$ and the remaining chord system of $\alpha$ is $\partial$-relatively 
homotopic to the remaining chord system of $\alpha'$ in $S^2$. 
Homotopic loop-chord diagrams are equivalent. 
The following definition is basic to this paper.

\phantom{x}

\phantom{x}

\noindent{\bf Definition~1.1.}  A connected loop-chord system $(o,\alpha)$ in $S^3$ is  
{\it circular}  if there are the same number of  loops  $o_i\,(i=1,2,\dots, n)$ of $o$ and 
chords $\alpha_i\,(i=1,2,\dots,n)$ of $\alpha$ such that the chord $\alpha_i$ connects $o_i$ to $o_{i+1}$ for every $i$  with $n+1\equiv 1$. 
A loop-chord system is {\it circular} if every 
connected component is a circular loop-chord system.
A disk-chord system $(d,\alpha)$ is {\it circular} if the loop-chord system $(o,\alpha)$ is 
circular, and {\it circular primitive} or briefly {\it CP} if every chord of $\alpha$ intersects
the interior of every disk in $d$ transversely meets $\alpha$ at just one point. 
 
\phantom{x}

\begin{figure}[hbtp]
\begin{center}
\includegraphics[width=12 cm, height=5 cm]{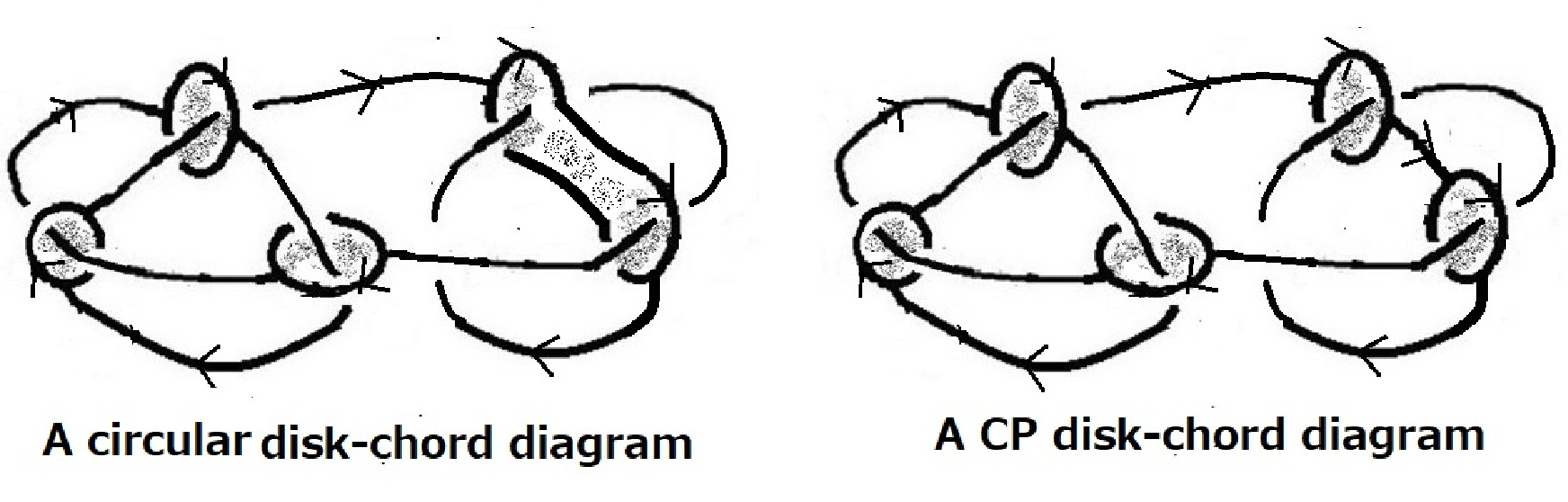}
\end{center}
\caption{An example of a circular disk-chord diagram and a CP disk-chord diagram that are equivalent by the move $M_1$}
\label{fig:CandCP}
\end{figure} 

\begin{figure}[hbtp]
\begin{center}
\includegraphics[width=6.5 cm, height=3 cm]{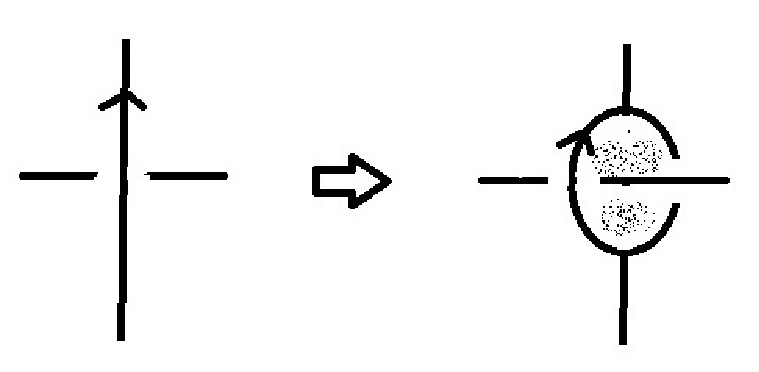}
\end{center}
\caption{Changing the location from near a crossing point to near 
a loop with chords or a disk with chords}
\label{fig:change}
\end{figure}

\begin{figure}[hbtp]
\begin{center}
\includegraphics[width=12cm, height=5 cm]{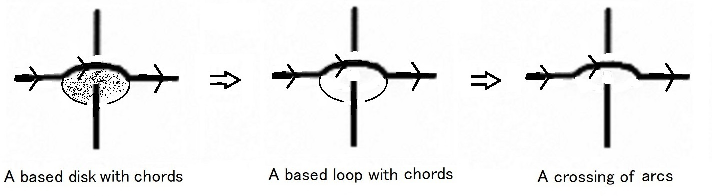}
\end{center}
\caption{A link $k$ with a disk system $d$}
\label{fig:disklink}
\end{figure}

A {\it CP disk-chord diagram} on $S^2$ is a diagram $C(d,\alpha)$ on $S^2$ 
of a CP disk-chord system $(d,\alpha)$ in $S^3$. 
Every ribbon torus-link $T$ in $S^4$ admits a CP disk-chord diagram  $C(d,\alpha)$, which is shown by an equivalent deformation of any given disk-chord diagram $C(d',\alpha')$ of $T$ using a finite number of the moves $M_0$, $M_1$ and $M_2$. 
An example of a circular disk-chord diagram and a CP disk-chord diagram that are equivalent by the move $M_1$ is given in Fig.~\ref{fig:CandCP}.
A loop-chord system $(o,\alpha)$ is also  {\it CP} if  the disk-chord system $(d,\alpha)$ is CP. The chord systems $\alpha$ of a circular loop-chord diagram $C(o,\alpha)$ and a circular loop-chord diagram $C(d,\alpha)$ have the unique orientation induced from the clockwise orientation of $o$.                                                                
From a disk-chord diagram $C(d,\alpha)$, the sphere-chord system $(O,\alpha)$ in $S^4$ is uniquely constructed by first realizing as a disk-chord system $(d,\alpha)$ in $S^3$ and 
then taking the trivial $S^2$-link $O=\partial(d\times[-1,1])$ for a collar $d\times[-1,1]$ 
of the disk system $d$ in $S^4$ given by  a collar $S^3\times[-1,1]$ of $S^3$ in $S^4$. 
From this viewpoint, the sphere-chord system $(O,\alpha)$ is {\it circular} if 
if the disk-chord system $(d,\alpha)$ is circular, and
the statement \lq\lq the chord system $\alpha$ meets transversely the interior of a disk of $d$ at  just $n$ points\rq\rq is expressed as \lq\lq the corresponding sphere component of the based sphere system $O$ is passed through $n$ times by  the chord 
system $\alpha$\rq\rq.
The spun torus-link $T(k)$ in $S^4$ for every link $k$ in $S^3$ is equivalent to 
the ribbon torus-link $F(CD)$ of  the unique CP disk-chord diagram $CD=C(d,\alpha)$ 
transformed from any given diagram $D$ of the link $k$ by the 
transformation given in Fig.~\ref{fig:change}.  This is explained  in  \cite{[4]}, and  for convenience in Section~2. Note that the diagram $D$ is recovered from 
the CP disk-chord diagram $CD$ by taking the upper arc system of the based loop system 
$o$, Fig.~\ref{fig:disklink}. The following theorem is a main result. 

\phantom{x}

\noindent{\bf Theorem~1.2.}  Every CP disk-chord diagram $C(d,\alpha)$ on $S^2$ of 
the spun torus-link $T(k)$ in $S^4$ of a link $k$ in $S^3$ is homotopic to the disk-chord diagram $CD$ of a diagram $D$ on $S^2$ of  $k$.  

\phantom{x}

This proof of Theorem~1.2 is useful to clarify  the method used in \cite{[8]}.
A double point of a diagram $D$ is called a {\it crossing point} of $D$. 
The number of crossing points in  $D$ is called the {\it crossing number} of $D$ 
and denoted by $c(D)$. 
The {\it crossing number} $c(k)$ of a link $k$ is 
the minimum of the crossing numbers $c(D')$ of  all diagrams $D'$ of all links $k'$ equivalent the link $k$. 
The {\it chord index} of a circular disk-chord diagram $C(d,\alpha)$ is 
the number $I(d,\alpha)$ of 
the geometric transverse intersection number of the interior of the disk system $d$ and 
the chord system $\alpha$ of the disk-chord system $(d,\alpha)$ in $S^3$.  
The {\it chord index} of a ribbon torus-link $F$ is the minimum $I(F)$ of the chord indexes  $I(d',\alpha')$ of all circular disk-chord diagrams $C(d',\alpha')$ of ribbon torus-links $F'$ equivalent to the ribbon torus-link  $F$.
If the chord index $I(F)$ of a ribbon torus-link $F$ is $n$, then there is a CP  
disk-chord diagram $(d,\alpha)$ of chord index $n$ of $F$ because a circular 
disk-chord diagram of chord index $n$ is equivalent to a CP disk-chord diagram of 
chord index $n$ by Fusion-fission move $M_1$. 
Since the span torus-link $T(k)$ has the CP disk-chord diagram $CD$ 
of any diagram $D$ of $k$, the inequality $c(k)\geq I(T(k))$ holds. 
By Theorem~1.1, $c(k)\leq I(T(k))$. 
Thus, the following corollary is obtained. 

\phantom{x}

\noindent{\bf Corollary~1.3.} $c(k)=I(T(k))$ for every link $k$ and the spun torus-link $T(k)$.

\phantom{x}

The additivity of crossing numbers  is known for  large classes of knots, e.g., alternating knots, adequate knots, torus knots, and investigated by lots of researchers, \cite{[9],  [10],  [11], [12], [13], [14]}.  By further developing the argument on Theorem~1.1, the following general additivity of the crossing numbers can be shown.

\phantom{x}

\noindent{\bf Theorem~1.4.} Let $k_1\# k_2$ be a connected sum of  
any knots or links $k_i\, (i=1,2)$. 
Then the identity $c(k_1\# k_2)=c(k_1)+c(k_2)$ holds.

\phantom{x}

\noindent{\bf 2. Proof of Theorem~1.2.}

Let $B$ be a 3-ball in $S^3$ with boundary sphere  $S=\partial B$.
The spun torus-link of a link $k$ in $B$ is the torus-link 
$T(k)=k\times S^1$ in the 4-sphere $S^4=B\times S^1\cup S\times D^2$, 
where $S^1$ and $D^2$ denote the unit circle and the unit disk in 
the complex number plane, respectively. 
For a boundary collar $S\times[0,1]$ of $S$ in $B$ with $S\times\{0\}=S$,  
assume that the link $k$ is in $S\times[0,1]$ and  the image of $k$ under 
the projection $S\times[0,1]\to S$ onto the first factor $S$ is a diagram $D$ of $k$. 
Write every point of $k$ as a pair $(x, s)$ with $x\in S$ and $s\in[0,1]$.
Then the {\it projection trace} of $k$ is the immersed annulus system (namely, the union of transversely meeting annuli) 
$P(k)=\cup_{(x,s)\in k} \{x\}\times[0,s]$ 
with $P(k)\cap S=D$. 
The union $V=P(k)\times S^1\cup D\times D^2$ is an immersed solid torus system 
in $S^4$ with ribbon disk singularity whose boundary is $T(k)$.   A multi-punctured 
solid torus system $W$  is obtained from $V$ by removing an open 3-ball 
neighborhood of the interior disk of every ribbon disk singularity, so that 
$\partial W=T(k)\cup O$ for a trivial $S^2$-link $O$ in $S^4$.  
The multi-punctured solid torus system $W$ is uniquely determined by 
a CP disk-chord system $(d,\alpha)$ in $S\times[0,1]$ with 
the intersection arc system $k_d=(\partial d)\cap k$  and the complementary arc system 
$\alpha=\mbox{cl}(k\setminus k_d)$.      
The multi-punctured solid torus system $W$ is directly constructed from the CP 
handled sphere-chord system $(O,h)$ with $h$ a 1-handle system with core system $\alpha$ obtained from the disk-chord system $D(d,\alpha)$ of a diagram of $k$  so that  $W=O\times[0,1]\cup h$ in $S^4$ 
with $\partial W=T(k)\cup O\times\{1\}$, where $O\times\{0\}=O$. 
From this viewpoint, the multi-punctured solid torus system $W$ is called 
a {\it SUPH system} for the ribbon torus-link $T(k)$ in $S^4$, \cite{[15]}. 
Let $\delta$ be a disk in $S$ with $\delta^c=\mbox{cl}(S\setminus\delta)$ 
the complementary disk of $\delta$ in $S$. 
Assume that the intersection $t=k\cap(\delta\times[0,1])$ is a {\it tangle} in the 3-ball 
$\delta\times[0,1]$, namely a proper arc system in  $\delta\times[0,1]$ 
such that the boundary point set $\partial t$ in the circle 
$(\partial \delta)\times\{0\}=\partial \delta$. 
Then the projection  $\delta\times[0,1]\to \delta$ sends the tangle $t$ to 
a tangle diagram $t^D$ in $\delta$ with $\partial t^D=\partial t$, which is a subdiagram of 
the diagram $D$.  
The {\it spun \mbox{$S^2$}-link}  of a tangle $t$ in $\delta\times[0,1]$ is 
the $S^2$-link$\mbox{cl}T(t)$ in $S^4$ given by the union 
$T(t)\cup (\partial t)\times D^2$. 
For the {\it projection trace} of the tange $t$, that is the union 
$P(t)=\cup_{(x,s)\in t} \{x\}\times[0,s]$, 
the union $V(t)=P(t)\times S^1\cup D\times D^2$ is an immersed 3-ball system with 
ribbon disk singularity whose boundary is $\mbox{cl}T(t)$. 
 A multi-punctured 3-ball system $W(t)$
of  the CP disk-chord diagram $D(d_t,\alpha_t)$ of 
the diagram $t^D$ of $t$ is obtained from $V(k)$ by removing an open 3-ball neighborhood 
of the interior disk of every ribbon disk singularity.  
The {\it cross-index} of a tangle diagram $t^D$ is the total number $\varepsilon(t^D)$ 
of  distinct arc pairs in $t^D$ such that each arc pair has 
the $Z_2$-intersection number $1$ in $\delta$, where   
define $\varepsilon(t^D)=0$ for a single arc diagram $t^D$, \cite{[16]}. 
The following lemma is  used for the proof of Theorem~1.2.

\phantom{x}

\noindent{\bf Lemma~2.1.} 
For every tangle $t$ in $\delta\times[0,1]$ with  tangle diagram $t^D$ such that 
the cross-index $\varepsilon(t^D)$ is $n>0$, then  every CP disk-chord diagram 
$t^D(d_t,\alpha)$ of the  spun $S^2$-link $\mbox{cl}T(t)$  must have at least 
$n$ based spheres.

\phantom{x}

\noindent{\bf Proof.} Let $(t_1, t^2)$ be a pair of distinct strings of the tangle $t$ 
such that the diagram pair $(t_{1p}^D,t_{2p}^D)$ in $\delta$ has 
the $Z_2$-intersection number $1$. 
For every transverse intersection point $p$ between $t_1^D$ and $t_2^D$, 
take the point $(p,t_{ip})$  in  $t_i\,(i=1,2)$. Say $t_{1p}<t_{2p}$. Then $t_{1p}$ belongs 
to the half-open interval $(t_{2p},0]$. Then the arc $t_1$ intersects the immersed 3-ball $V(t_2)$ transversely at the point $(p,t_{1p})$. By counting the case $t_{1p}>t_{2p}$, 
it is concluded that  the $Z_2$-intersection number  of $(t_{1p}^D,t_{2p}^D)$  is $1$ 
in $\delta$ if and only if 
the sum of  the $Z_2$-intersection numbers of $(t_1, V(t_2))$ and $(t_2, V(t_1))$
is $1$. This claim holds regardless of choices of immersed 3-balls $V(t_i)$ 
bounded by $\mbox{cl}T(t_i)\,(i=1,2)$ since the $Z_2$-linking numbers of 
$(t_1, \mbox{cl}T(t_2))$ and $(t_2, \mbox{cl}T(t_1))$ are well-defined. 
Thus, for $\varepsilon(t^D)=n$, every CP disk-chord diagram 
$Ct^D$ of the  spun $S^2$-link $\mbox{cl}T(t)$  must have at least 
$n$ based spheres. This completes the proof of Lemma~2.1.

\phantom{x}

The proof of Theorem~1.2 is done as follows.

\phantom{x}

\noindent{\bf 2.2: Proof of Theorem~1.2.} 
Let  $C(d,\alpha)$ be the CP disk-chord diagram on the sphere $S=\partial B$
of the spun torus-link $T(k)$ in $S^4$, and $(d,\alpha)$ a CP disk-chord system 
in the interior of  a boundary-collar $S\times[0,1]$ of $S$ in $B$ which is sent to  
the CP disk-chord diagram$C(d,\alpha)$ on $S$ under the projection 
$S\times[0,1]\to S$. 
Assume that the union $d\cup\alpha$ is connected since otherwise the task can be done 
for each connected component.
Further, by assuming that the spun torus-link $T(k)$ is constructed from a link $k$ in the interior of $S\times[0,1]$ and then deforming the union 
$d\cup\alpha$ in $S^4$ into $B$ by an ambient isotopy of $S^4$ keeping the link $k=k\times\{1\}\subset T(k)$ fixed, the union 
$d\cup\alpha$ is written as the union $d\cup k$ in $S\times[0,1]$ by counting the relationship between $d\cup\alpha$ and $d\cup k$ around the based disk system $d$, illustrated in Fig.~\ref{fig:disklink} in $S\times[0,1]$.  
Let $\Gamma$ be a connected graph in $S\times [0,1]$ 
obtained from $d\cup k$ by shrinking every disk 
of $d$ into a vertex. 
Let $(\Delta, \tau^+)$ be a regular neighborhood of $\tau$ in 
$(S\times[0,1],\tau)$ with $\Delta$ a 3-ball.
Let $e=\mbox{cl}(\Gamma\setminus\tau^+)$ and $t=\mbox{cl}(k\setminus e)$ be 
be arc systems in $k$. Let 
 $\tau^+ (d,t)$ be the disk-arc system obtained by replacing 
the vertex system of $\Gamma$ with the disk system $d$. 
Deform the 3-ball $\Delta$ into the 3-ball $\delta\times[0,1]$ in $S\times[0,1]$ 
for a disk $\delta$ in $S$, so that the projection $\delta\times[0,1]to\delta$ sends 
the disk-tangle system $\tau^+ (d,t)$ into  the disk-tangle diagram $D\tau^+ (d,t)$ 
in $\delta$ and the boundary point system $\partial t$ of $t$ is in the boundary circle 
$\partial\delta$ of $\delta$. Note that the disk-tangle system $D\tau^+ (d,t)$ is 
a disk-chord diagram of a tangle diagram $t^D$ in $\delta$. 
The arc system $e$ is a tangle in the 3-ball 
$\delta^c\times[0,1]$ for the disk $\delta^c=\mbox{cl}(S\setminus\delta)$ 
with the boundary point system $\partial e$ in the boundary circle $\partial \delta^c$. 
By construction, the spun $S^2$-link $\mbox{cl}T(t)$ bounds 
an immersed 3-ball system $V(t)$ with ribbon disk singularity  given by  the disk-tangle diagram $D\tau^+ (d,t)$ and 
the spun $S^2$-link $\mbox{cl}T(e)$ bounds a solid torus system $V(e)$ such that 
the union $V(t)\cup V(e)$ is an immersed handlebody system constructed from
the CP disk-chord diagram $C(d,\alpha)$ on the sphere $S$.
By Lemma~2.1, the cross-index $\varepsilon(e^D)$ of the diagram $e^D$ of $e$ in 
the disk $\delta^c$ is $0$. 
Since the spun $S^2$-link $\mbox{cl}T(e)$ of the tangle $e$ 
in the 3-ball $\delta^c\times[0,1]$ is a trivial $S^2$-link in $S^4$, 
the tangle $e$ is a trivial tangle in $\delta^c\times[0,1]$. 
In fact, the compact exterior $E$ of $e$ in  $\delta^c\times[0,1]$ has a meridian-based 
free fundamental group $\pi_1(E,b_0)$ and $E$ is a handlebody by Dehn's lemma, \cite{[1]}.
Thus, there is a disjoint 
disk system $d_e$ in $\delta^c\times[0,1]$ with $\partial d_e=e\cup e^c$ 
for a disjoint arc system $e^c$ in the annulus $(\partial\delta^c) \times[0,1]$.
The cross-index $\varepsilon(e^D)=0$ implies that the arc system $e$ is isotopic into
an arc system in $\delta^c$ by an isotopy of $\delta^c\times[0,1]$ keeping the boundary point system $\partial e$ fixed. Thus, there is a link diagram $D$ of $k$ on $S$ such that 
the CP disk-chord diagram $C(d,\alpha)$ is the disk-chord diagram 
$CD$ of a link diagram $D$. 
This completes the proof of Theorem~1.2.

\phantom{x}

\noindent{\bf 3. Proof of Theorem~1.4.} 

If $k_1$  is a trivial knot, then the identity $c(k_1\# k_2)=c(k_1)+c(k_2)$ holds. 
If $k_1$ is a split link of sublinks $k'_1$ and $k''_1$ and a connected sum $k_1\# k_2$ is 
made as $k'_1\# k_2\cup k''_1$ and the identity $c(k'_1\# k_2)=c(k'_1)+c(k_2)$ is known, then 
$c(k_1\# k_2)=c(k_1)+c(k_2)$ also holds. 
Thus, for the proof of Theorem~1.4, it may be assumed that the links $k_i\, (i=1,2)$ are non-trivial and non-split links. Then every component of every connected sum $k_1\# k_2$ does not bound any singular disk with interior disjoint from $k_1\# k_2$, by Dehn's lemma, \cite{[1], [17]}. 
Under this assumption, the ribbon torus-link $T(k_1\# k_2)$ admits a {\it unique meridian system} in $S^4$, that is, there is a unique simple loop up to homotopy in every component of the ribbon torus-link 
$T(k_1\# k_2)$  which bounds  a  disk  in $S^4$ with 
the interior disjoint from $T(k_1\# k_2)$, 
because there is a natural isomorphism $\pi_1(S^3\setminus k_1\# k_2, x_0)\to \pi_1(S^4\setminus T(k_1\# k_2), x_0)$. 
Let $D_i$ be a diagram of  a link $k_i$ with $c(D_i)=c(k_i)\,(i=1,2)$, 
and $D_{12}$ a diagram of  a connected sum $k_1\# k_2$ with $c(D_{12})=c(k_1\# k_2)$. 
Let $A_i\, (i=1,2)$ be 4-balls in $S^4$ with $A_1\cap A_2=B_0$ a 3-ball in the boundary 
3-spheres $\partial A_i\, (i=1,2)$. 
Take a handled-sphere system $(O_i, h_i)$ for the torus-link $T(k_i)$ in the 4-ball $A_i$  constructed from the disk-chord system $CD_i=C(d_i,\alpha_i)$ of the diagram $D_i$, 
where $h_i$ denotes a 1-handle system with core system $\alpha_i$ an arc system in $k_i$.  
Let $(O_{12}, h_{12})$ be a handled-sphere system for the torus-link $T(k_1\# k_2)$ in 
the 4-ball $A_1\cup A_2$ 
constructed from the CP disk-chord system $CD_{12}$ of the diagram $D_{12}$, 
where $h_{12}$ denotes a 1-handle system with core system $\alpha_{12}$ an arc system in $k_1\# k_2$.
By deforming $\alpha_i$ in $A_i\,(i=1,2)$, assume that the chord systems 
$\alpha_i\,(i=1,2)$ meets  the 3-ball $B_0$ with an arc $\|a_0 b_0\|$ with endpoint pair 
$(a_0,b_0)$.   
Let $O'_{12}=O_1\cup O_2$ and 
$\alpha'_{12}=\mbox{cl}(\alpha_1\setminus \|a_0b_0\|)\cup \mbox{cl}(\alpha_2\setminus 
\|a_0b_0\|)$.  
Then there is a handled-sphere system $(O'_{12}, h'_{12})$ for $T(k_1\#k_2)$ where 
$h'_{12}$ denotes a 1-handle system with core system $\alpha'_{12}$. 
meeting the 3-ball $B_0$ with two transverse disks $\mu, \nu$ of a 1-handle 
in $h'_{12}$. 
Under these preliminaries, the proof of Theorem~1.4 is done as follows.

\phantom{x}

\noindent{\bf Proof of Theorem~1.4.}  
Consider that the connected sum $k_1\#k_2$ is done between a component $k^0_1$ of 
$k_1$ and a component $k^0_2$ of  $k_2$ and let $k'_i=k_i\setminus k^0_i\,(i=1,2)$.
Let $W_{12}=O_{12}\times[0,1]\cup h_{12}$ and $W'_{12}=O'_{12}\times[0,1]\cup h'_{12}$ 
be the SUPH systems for $T(k_1\# k_2)$. 
After an isotopic  deformation of the chord system $h'_{12}$ in $A_1\cup A_2$ fixing 
$\|a_0b_0\|\cup O'_{12}$, there is a multi-punctured solid torus system $W'_*$ in 
$A_1\cup A_2$ 
obtained from $W'_{12}$ by adding a standard 2-handle system on $O'_{12}\times\{1\}$ 
such that there is an orientation-preserving diffeomorphism of $A_1\cup A_2$ 
sending  $W_{12}$ is sent to the manifold $W'_{**}$ obtained from $W'_{*}$ 
by removing an embedded 1-handle system $h''$ on the trivial $S^2$-link $O'_*=\partial W'_*\setminus T(k_1\#k_2)$, \cite[Appendix]{[15]}. 
Let $\delta''$ be the transverse disk system of the 1-handle system $h''$ cut by the
disk union $\mu\cup\nu$.
Assume $W'_{**}=W_{12}$. 
Let $(W_{12})_0$ be the SUPH system component  of $W_{12}$ containing the disk 
union $\mu\cup\nu$, and $(O_{12}, \alpha_{12})_0=((O_{12})_0, (\alpha_{12})_0)$  
the circular spherical chord system component of the CP spherical chord system 
$(O_{12}, \alpha_{12})$ in $(W_{12})_0$. 
In terms of $(O_{12}, \alpha_{12})_0$, the disk system $\delta''$ is a proper disk system 
in a 3-ball system $B_{12}$ in the 4-ball $A_1\cup A_2$ 
bounded by the trivial $S^2$-link $(O_{12})_0$ after 
the chord system $(\alpha_{12})_0$ is deformed by a $\partial$-relative isotopy in 
$A_1\cup A_2$, \cite[Proposition A.1]{[15]}.  
Since the circular disk-chord system $CD_{12}$ is  primitive,  
the intersection  $\delta''\cap B_{12}$ is changed into the intersection 
$\delta''\cap (h_{12})_0$ for the 1-handle system $(h_{12})_0$ with $(\alpha_{12})_0$ 
the core arc system. Thus, the intersection of  the handled sphere system 
$(O_{12}, h_{12})_0$ and the disk union $\mu\cup\nu$ can be assumed to be 
the intersection of $(h_{12})_0$ and  $\mu\cup\nu$.   
Then, since the compact exterior $E_{12}$ of $O_{12}\cup\alpha_{12}$ in $W_{12}$ is diffeomorphic to $T(k_1\# k_2)\times[0,1]$, 
$T(k_1\# k_2)$ has a unique meridian system in $S^4$, 
and $(O_{12}, \alpha_{12})_0$ is a circular spherical chord system of $(W_{12})_0$,  
the oriented chord system $(\alpha_{12})_0$ meets $\mu$ and $\nu$
with intersection numbers $+1$ and  $-1$  in $(W_{12})_0$, respectively. 
There is an band system  with a core arc system in  $\mu\cup\nu$ in $(W_{12})_0$
which spans chords in $(\alpha_{12})_0$ and  splits the connected circular sphere-chord system  $(O_{12}, \alpha_{12})_0$
into a connected circular sphere-chord system $C_0$ meeting $\mu$ and $\nu$ at one point each and  some circular sphere-chord systems $C_j\, (j=1,2,\dots,s)$ not meeting 
$\mu\cup\nu$, where $C_0$ is slided on $\mu$ and $\nu$ to set $C_0\cap \mu=a_0$ and 
$C_0\cap \nu=b_0$.   
Let $C'_0$ be a circular sphere-chord system in $(W_{12})_0$ by joining $C_0$ and 
$C_j\, (j=1,2,\dots,s)$ by a band not meeting $\mu\cup\nu$. 
Let  $(O^C_{12}, \alpha^C_{12})$ be a CP sphere-chord system in $W_{12}$ 
obtained from $(O_{12}, \alpha_{12})$ by replacing $(O_{12}, \alpha_{12})_0$ with 
$C'_0$, which is a circular sphere-chord system of $W_{12}$,  
because the 3-manifolds obtained from $(W_{12})_0$ by splitting along 
the disk system are simply connected. 
Let $W'_i\,(i=1,2)$ be the 3-manifolds obtained from $W_{12}$ by splitting along 
$\mu\cup\nu$, where a boundary collar of  $\mu\cup\nu$ in $W'_i$ belongs to $A_i$. 
Since $W'_i\,(i=1,2)$ are disjoint in $A_1\cup A_2$ but in general not split in $S^4$,  
move $W'_i$ into $A_i$ by disregarding the other 3-manifold $W'_{i'}\, (i'\ne i)$ 
and keeping $\mu\cup\nu$ fixed . 
After  this move, while $W_{12}$ remains a  SUPH system  for $T(k_1\# k_2)$ 
and $(O^C_{12}, \alpha^C_{12})$ is a circular sphere-chord system,  
but $(O^C_{12}, \alpha^C_{12})$ is no longer any CP sphere-chord system. 
Let $W'_i$ be a SUPH system for $T(k_i)$ obtained from $W'_i$ by adding a 
1-handle $h_0$ on $\mu\cup\nu$ with core arc $\|a_0 b_0\| \,(i=1,2)$, which is 
a SUPH system of  a circular sphere-chord system $(O^C_i, \alpha^C_i)$ 
obtained from $(O^C_{12},\alpha^C_{12})$ by splitting along a band with core arc 
$\| a_0 b_0\|\,(i=1,2)$. 
\[c(k_1\# k_2)=r(O_{12})=r(O^C_{12})=r(O^C_1)+r(O^C_2).\]
Every sphere of $O^C_i\,(i=1,2)$ is  passed through at most once time by the chord system $\alpha^C_i$ by construction. 
By Corollary~1.3, $r(O^C_i)\geq I(T(k_i))=c(k_i)\,(i=1,2)$. 
so that $c(k_1\# k_2)\geq c(k_1)+c(k_2)$. Since $c(k_1\# k_2)\leq c(k_1)+c(k_2)$ by definition, the identity  $c(k_1\# k_2)=c(k_1)+c(k_2)$ holds.
This completes the proof of Theorem~1.4.

\phantom{x}

\noindent{\bf Acknowledgments.} 
This work was partly supported by JSPS KAKENHI Grant Numbers JP21H00978
and JP26K06456, and MEXT Promotion of Distinctive Joint Research Center 
Program JPMXP0723833165 and Osaka Metropolitan University Strategic Research 
Promotion Project (Development of International Research Hubs).

\phantom{x}

,

\end{document}